\documentclass[11pt]{article}

\usepackage[a4paper, margin=1in]{geometry}
\usepackage{amsfonts,amsmath,amssymb,amsthm,graphicx,fixmath,latexsym,color}
\usepackage{xcolor}
\usepackage{hyperref,epsfig}
\usepackage{algorithmic}
\usepackage{algorithm}
\usepackage{xspace}

\usepackage{graphicx}
\usepackage{caption,subcaption}

\providecommand{\remove}[1]{}

\newtheorem{theorem}{Theorem}[section]

\newtheorem{claim}[theorem]{Claim}

\newtheorem{conjecture}[theorem]{Conjecture}

\newtheorem*{theorem*}{Theorem}
\newtheorem*{lemma*}{Lemma}
\newtheorem*{proposition*}{Proposition}

\newtheorem*{example*}{Example}

\newcommand{\bbox}{\vrule height7pt width4pt depth1pt}

\begin{document}
\title{On sets of $n$ points in general position that determine
lines that can be pierced by $n$ points}

\author{Chaya Keller\thanks{Mathematics Department, Technion -- Israel Institute of Technology, Haifa 32000, Israel. \texttt{chayak@technion.ac.il}. Research partially
supported by Grant 409/16 from the Israel Science Foundation.}
\and
Rom Pinchasi\thanks{Mathematics Department, Technion -- Israel Institute of Technology, Haifa 32000, Israel. \texttt{room@tx.technion.ac.il}. Research partially supported by Grant 409/16 from the Israel Science Foundation.}
}

\date{}
\maketitle

\begin{abstract}
Let $P$ be a set of $n$ points in general position in the plane.
Let $R$ be a set of $n$ points disjoint from $P$ such that
for every $x,y \in P$ the line through $x$ and $y$ contains a point in $R$
outside of the segment delimited by $x$ and $y$. We show that
$P \cup R$ must be contained in a cubic curve.
\end{abstract}


\section{Introduction}
\label{sec:introduction}

A beautiful result of Motzkin \cite{M67}, Rabin, and Chakerian \cite{Ch70}
states that
any set of non-collinear red and blue points in the plane determines a
monochromatic line. Gr\"unbaum and Motzkin \cite{Gr75} initiated the study
of biased coloring, that is, coloring of the points such that no purely
blue line
is determined. The intuition behind this study is that if the number of blue
points is much larger than the number of red points, then unless the set
of blue points is collinear, the set of blue and red points should determine a
monochromatic blue line.

The same problem was independently considered by Erd\H os and Purdy
\cite{EP78} who stated
it in a slightly different way.

Let $P$ be a set of $n$ points in the plane. Erd\H os and Purdy
asked the following question in \cite{EP78}:
Assume a set $R$ of points in the plane is disjoint from $P$ and
has the property that every line
determined by $P$ passes through a point in $R$. How small can be the
cardinality of $R$ in terms of $n$?

Clearly, if $P$ is contained in a line, then $R$ may consist of just one
point. Therefore, the question of Erd\H os and Purdy is about sets $P$
that are not collinear.
The best known lower bound for this question is given
in \cite{Pinchasi13}, where it is shown that $|R| \geq n/3$.

\begin{figure}[ht]
    \centering
    \includegraphics[width=7cm]{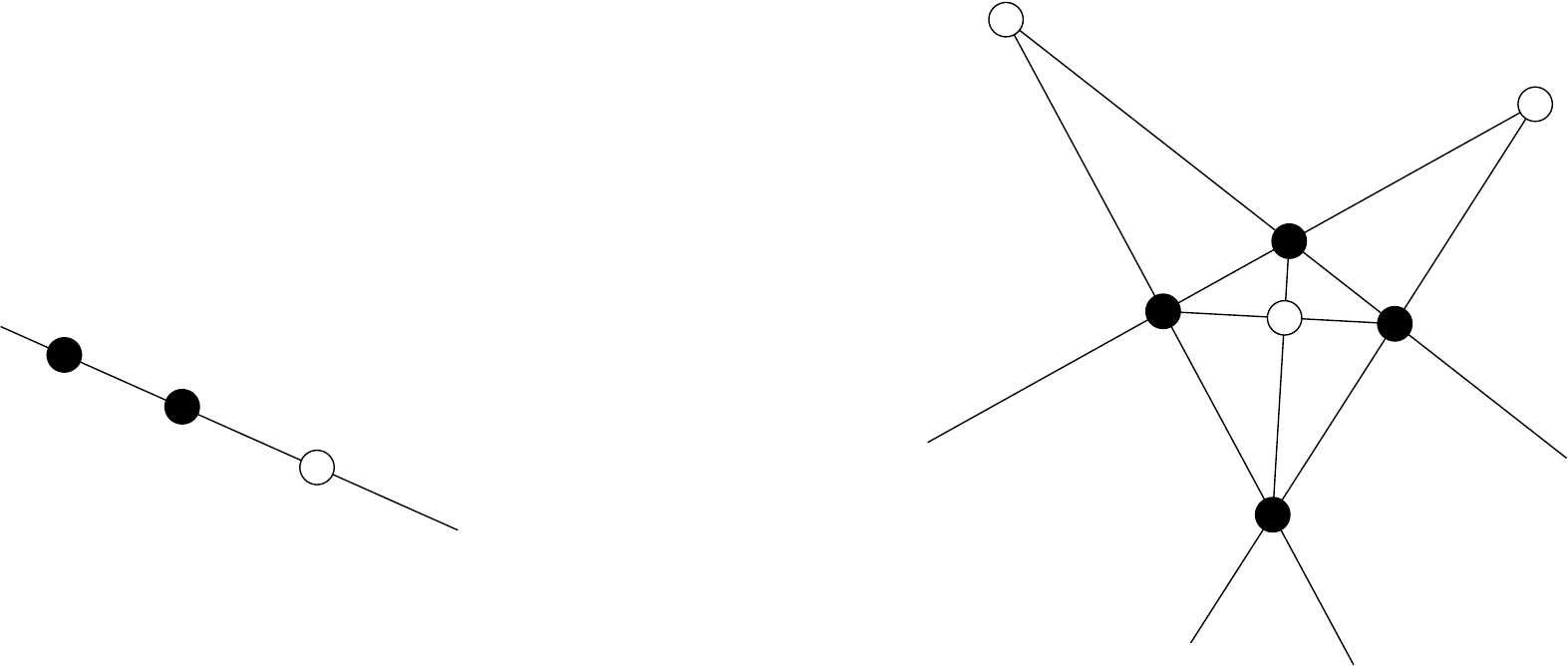}
        \caption{Constructions with $|R|=n-1$ for $n=2,4$.
The points in $P$ are colored black while the points in $R$ are colored white.}
        \label{figure:fig1}
\end{figure}

In \cite{EP78} Erd\H os and Purdy considered the same problem
in the case where the set $P$
is in \emph{general position} in the sense that no three points of $P$
are collinear. In this case if $n$ is odd the tight bound $|R| \geq n$
is almost trivial because every point in $R$ may be incident to at most
$\frac{n-1}{2}$ of the ${n \choose 2}$ lines determined by $P$.
To observe that this bound is tight let $P$ be the
set of vertices
of a regular $n$-gon and let $R$ be the set of $n$ points on the line
at infinity that correspond to the directions of the edges (and diagonals)
of $P$. This construction is valid also when $n$ is even.

If $n$ is even, a trivial counting argument shows that $|R|$ must be at least
$n-1$. This is because every point in $R$ may be incident to at most
$n/2$ lines determined by $P$. This trivial lower bound for $|R|$ is in fact
sharp in the cases $n=2$ and $n=4$, as can be seen in Figure \ref{figure:fig1}.
Is the bound $|R| \geq n-1$ sharp also for larger values of $n$?

The following theorem proves a conjecture attributed to Erd\H os and Purdy \cite{EP78}. This conjecture was addressed in many papers
(see, e.g.,~\cite[Chapter~7.3]{BMP} and the references therein).

\begin{theorem}[\cite{ABKPR08,M18,P18,PP19}]\label{conjecture:EP}
Let $P$ is a set of $n$ points in general position in the plane
and, where $n>4$ is even. Assume $R$ is another set of points disjoint from
$P$ such that every line through two points of $P$ contains a point from $R$.
Then $|R| \geq n$.
\end{theorem}

Theorem \ref{conjecture:EP} was first proved in
\cite{ABKPR08} (see Theorem 8 there), as a special case of the solution of the
Magic Configurations conjecture of Murty \cite{Murty71}. The proof
in \cite{ABKPR08} contains a topological argument based on Euler's formula
for planar maps and the discharging method.
An elementary (and long) proof of Theorem
\ref{conjecture:EP} was given by Mili\'{c}evi\'{c} in \cite{M18}.
An algebraic proof of Theorem \ref{conjecture:EP} is given
in \cite{P18}.
Probably the ``book proof'' of the Theorem \ref{conjecture:EP} can be found in
\cite{PP19}.

Theorem \ref{conjecture:EP} was proved also over $\mathbb{F}_{p}$
by Blokhuis, Marino, and Mazzocca \cite{BMM14}.

As we have seen, there are constructions of sets $P$ of $n$ points in general
position and sets $R$ of $n$ points not in $P$, such that every line
determined by $P$ passes through a point in $R$. One major question that
arises here is to characterize those sets $P$ in general position for which
there exists a set $R$ with $|R|=|P|$ such that every line that is determined
by $P$ passes through a point in $R$.

Already in \cite{EP78} Erd\H os and Purdy drew a connection between
the problem of finding a small set $R$ that pierces all the lines determined
by a set $P$ and another well-known problem, raised by Scott~\cite{Scott70} and called `the slopes problem':

\medskip \emph{What is the minimal number of directions spanned by a set $P$ of $n$ non-collinear points in the plane}?

This is equivalent (by a projective transformation) to finding a set $R$ of minimum size that is
\emph{contained
in a line}, and pierces all lines determined by a set $P$ of $n$ non-collinear
points.

The slopes problem was solved by Ungar~\cite{Ungar82}, in a beautiful proof.
Ungar showed that the minimum number of distinct direction of lines
determined by a set of $n$ non-collinear points in the plane is
$2\lfloor n/2 \rfloor$, and the bound is tight, e.g., for the vertex set of a regular $n$-gon for $n$ even, and for the vertex set of a regular $(n-1)$-gon plus its center point, for $n$ odd.

The extremal configurations for the slopes problem, i.e., sets of $n$ non-collinear points that span exactly $n-1$ directions, were studied in a series of papers of Jamison and Hill (see~\cite{Jamison85} and the references therein). They found four infinite families of configurations and $102$ sporadic examples, and their list is conjectured to be complete for $n \geq 49$.

Elekes~\cite{Elekes99} studied `almost extremal' configurations, for which at most $cn$ directions are spanned. He provided a characterization of all such configurations with at least $c'n$ points lying on a line and at least $c'' n$ points lying outside that line. In the same paper, Elekes conjectured that for any $m \geq 6$ and $c>0$, there exists $n_0$ such that any set of $n \geq n_0$ non-collinear points in the plane that determines at most $cn$ directions, contains at least $m$ points that lie on a quadratic curve. This conjecture is still wide open, even for $m = 6$.

When $P$ is assumed to be in general position in the sense that no three
of its points are collinear, then the slopes problem becomes almost trivial.
Every set of $n$ points in general position in the plane determines at least $n$
distinct directions for $n \geq 3$.

Extremal and almost-extremal configurations for the slopes problem
were studied also in the case where the point set in question is in general
position. Jamison~\cite{Jamison86} provided a complete characterization of the
extremal configurations in this case, showing that if $P$ is a set of
$n$ points in general position in the plane that determines exactly $n$
distinct directions of spanned lines then $P$ is an affine image of
vertex set of a regular $n$-gon. In particular $P$ is contained in a quadric and
because in this case the set $R$ is collinear, then $P \cup R$ is contained
in a cubic curve in the plane. This simple observation will directly
be related to the main result in this paper.

Jamison conjectured that if a set $P$ of $n$ points in general position
in the plane determines precisely $2n-c$ distinct direction
where $c_{0}  \leq  c \leq n$ for an absolute constant $c_0$,
then $P$ is obtained, up to an affine transformation,
from a regular $(2n-c)$-gon by omitting $n-c$ of its vertices.
Recently, Pilatte~\cite{Pilatte18} confirmed this conjecture for the case where $P$ determines precisely $n+1$ distinct directions.

\bigskip

We conjecture the following

\begin{conjecture}\label{conjecture:main}
Let $P$ be a set of $n$ points in general position and and let $R$
be a set of $n$ points disjoint from $P$. If every line determined by
$P$ passes through a point in $R$, then $P \cup R$ is contained in
a cubic curve.
\end{conjecture}

It would be tempting to conjecture an even much more far reaching connection
between the structure of $P$ and the cardinality of $R$, namely that if
$|R|=O(|P|)$, then ``many'' points of $P$ lie on some low degree (cubic?)
algebraic curve. This would correspond to the conjecture of Elekes
mentioned above.

The main result in this paper is a proof of
Conjecture \ref{conjecture:main} under the additional assumption that the
line through every two points $x,y \in P$
contains a point from $R$ that lies outside of the segment
determined by $x$ and $y$.

\begin{theorem}\label{theorem:main}
Suppose $P$ is a set of $n$ points in general position in the plane and
$R$ is another set of $n$ points, disjoint from $P$.
If for every $x,y \in P$ there is a point $r \in R$ on the line through $x$
and $y$ and outside the interval determined by $x$ and $y$, then
$P \cup R$ is contained in a cubic curve.
\end{theorem}

Theorem \ref{theorem:main} generalizes in particular the characterization
of Jamison in \cite{Jamison86} for sets in general position that determine
precisely $n$ distinct directions. Indeed, suppose a set $P$ of $n \geq 4$ points
in general position in the plane determines precisely $n$ distinct directions,
then there exists a collinear set $R$ of $n$ points (on the line at infinity)
such that every line determined by two points $x$ and $y$ in $P$ contains a
point of $R$.
This point of $R$ must lie outside the segment delimited by $x$ and $y$
because $R$ is contained in the line at infinity.
By Theorem \ref{theorem:main} $P \cup R$ is contained in a cubic $C$. Because
$n \geq 4$, $C$ contains at least $4$ collinear points on the line at infinity
and therefore contains all the line at infinity. It follows now
that $C$ is a union of a quadric and the line at infinity. Therefore,
$P$ must be contained in a quadric.

The proof of Theorem \ref{theorem:main} has two steps. In the first step,
that is more of a combinatorial nature,
we use our assumption that the point of $R$ collinear with two
points $x$ and $y$ in $P$, lies outside of the segment delimited by $x$ and $y$.
We show that in such a case the points of $P$ must lie in convex position.
Moreover, the structure of collinearities in $P \cup R$ is uniquely determined:
if we denote the points of $P$ by $x_{0}, \ldots, x_{n-1}$
in the cyclic order in which they are arranged in convex position, then
one can rename the points in $R$ by $r_{0}, \ldots, r_{n-1}$ such that
two points $x_{i}$ and $x_{j}$ in $P$ are collinear with $r_{k}$ iff
$i+j+k=0$ modulo $n$.
In the second step of the proof we use this information and follow
the footsteps of Green and Tao in \cite{GT13}, where the main algebraic ingredient
is Chasles' theorem (see below), to conclude that $P \cup R$
is contained in some cubic curve.

\section{Proof of Theorem \ref{theorem:main}}

Let $P'$ denote the set of points in $P$ that are extreme on the convex hull
of $P$. Let $k=|P'|$.
We first show that there are precisely $k$ points of $R$ outside of the convex
hull of $P$.
Fix a point $x$ in $P$. For every $y \in P\setminus \{x\}$
there is a point in $r \in R$ collinear with $x$ and $y$ that lies
outside of the segment delimited by $x$ and $y$. In principle
there can be more than one such point $r$. We choose one and say that it is
\emph{relevant} for $x$. Observe that for every $x$ in $P$ there is precisely
one point in $R$ that is not relevant for $x$.

As a consequence, a line $\ell$ through two points in $P$ may contain only
one or two points of $R$. This is because if it contains three, then for each
of the two points of $P$ on $\ell$ there are at least two points of $R$ that
are not relevant for it, contrary to our observation.

For every two points $x,y \in P'$ it is true that the point in $R$ collinear
with $x$ and $y$ must lie outside of the convex hull of $P$. Therefore,
by Theorem \ref{conjecture:EP} the number of points of $R$ outside the
convex hull of $P$ must be at least $k$.
On the other hand, let $R'$ denote the set of all points of $R$
outside the convex hull of $P$. For every point $r \in R'$ there are two lines
through $r$ that are tangent to the convex hull of $P$.
Among all those lines, let $a$ denote the number of the lines that support the convex hull of
$P$ at a single vertex.
Such a line must pass through precisely one point of $R'$.
This is because any point of $R'$ on this line is not relevant for
the unique point of $P'$ on this line and for every point of $P$ there is
precisely one point of $R$ that is not relevant for it.
Notice that $a \leq 2|R'|-k$ because there are $k$
lines supporting the convex hull of $P$ at an edge and each such line
must contain at least one point of $R'$.
Let $b$ denote the number of lines that support the convex hull of $P$ at
an edge and contain only one point of $R$.
Let $c$ denote the number of lines that support the convex hull of $P$ at
an edge and contain precisely two points of $R$.
We have $a+b+2c=2|R'|$ and $b+c=k$.

It follows that the number of times a point in $R'$ is not relevant for
a point in $P'$ is equal to $a+2c$. This should be
at most the cardinality of $P'$, namely $k$.

Therefore,
$k \geq a+2c =2|R'|-b \geq 2|R'|-b-c=2|R'|-k$
This readily implies $|R'| \leq k$.

By showing that $|R'| \geq k$ and $|R'| \leq k$, we conclude that
$|R'|=k$. Therefore, $c$ must be equal to $0$.
Moreover, for every point $p' \in P'$ there is a point $r' \in R'$
that is not relevant for it and the line through $p'$ and $r'$ is tangent to
the convex hull of $P'$ and supports it at the point $p'$ alone.

\begin{figure}[tb]
\begin{center}
\scalebox{0.6}{
\includegraphics{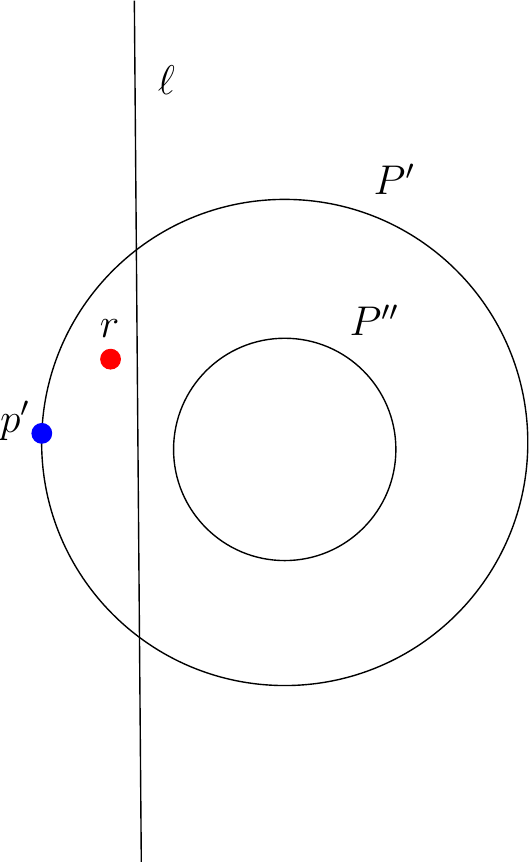}
}
\caption{Illustration for the proof that no point of $R$ lies inside the convex hull of $P'$ and outside the convex hull of $P''$.}
\label{fig:empty1}
\end{center}
\end{figure}

We will now show that $P'=P$. Assume not and let $P''$ be the set of
extreme vertices of $P \setminus P'$. We claim that no point of $R$ can lie
inside the convex hull of $P'$ and outside the convex hull of $P''$.
Indeed (see Figure~\ref{fig:empty1} for an illustration),
if $r \in R$ is such a point, consider a line $\ell$
separating $r$ and $P''$. It must be that there is a point $p' \in P'$
in the half-plane bounded by $\ell$ and not containing $P''$, for otherwise
$r$ cannot belong to the convex hull of $P'$. Notice that $r$ cannot be
relevant for $p'$. This is a contradiction because there is already a point
in $R'$ that is not relevant for $p'$.

Now take any point $p' \in P'$ and consider a line $\ell$ through $p'$ that
supports the convex hull of $P''$ at a point $p'' \in P''$.
The line $\ell$ must contain a point in $r \in R$. The point $r$ cannot
lie inside the convex hull of $P'$ because there is no point of $R$
inside the convex hull of $P'$ and outside the convex hull of $P''$.
Therefore, $r$ must be in $R'$. But then the line $\ell$ through $r$ and $p'$
contains only one point of $P'$ while it is not a tangent line of the convex
hull of $P'$.

We conclude that $k=n$, $P=P'$, and $R=R'$. That is, the points of $P$ are
in convex position and the points of $R$ lie outside the convex hull of $P$.
We will now show that the collinear triples of points two from $P$ and one from
$R$ are uniquely determined.

Denote by $x_{0}, x_{1}, \ldots, x_{n-1}$ the points of $P$ in their counterclockwise
cyclic order on the boundary of the convex hull of $P$.
In the following analysis the summation of indices is done modulo $n$.
We claim that one can rename the points in $R$ to be $r_{0}, \ldots, r_{n-1}$
in such a way that for every $0 \leq i \neq j <n$ and $0 \leq k < n$
the point $x_{i}, x_{j}$, and $r_{k}$ are collinear if and only if
$i+j+k=0$ (modulo $n$ of course).

Indeed, assume not, then one can find two distinct indices
$i$ and $j$ such that $x_{i}$ and $x_{j}$ are colllinear with $r_{k}$
for some $k$ while $x_{i+1}$ and $x_{j-1}$ are not collinear with $r_{k}$
and $i+1 \neq j-1$.
Assume without loss of generality that $x_{i}$ lies between $r_{k}$
and $x_{j}$ (the case where $x_{j}$ lies between $r_{k}$ and $x_{i}$ is
similar). Rotate the line through $r_{k}, x_{i}$, and $x_{j}$ about
$r_{k}$ in the clockwise direction until the first time it meets a point
in $P$. This point will be either $x_{i+1}$ or $x_{j-1}$ but not both at once.
Then we find a line $\ell$ through $r_{k}$ and only one point of $P$ that
does not support the convex hull of $P$. This is impossible because
as we have seen every point $p \in P$ has a point in $r \in R$ that is not
relevant for $p$ such that the line through $r$ and $p$ supports the convex hull
of $P$ and this point $r$ is the unique point in $R$ that is not relevant for
$p$. This completes the first step of the proof of Theorem \ref{theorem:main}.

\subsection{Fitting $P \cup R$ into a cubic curve}

In the second step of the proof of Theorem \ref{theorem:main} we will show
that $P \cup R$ is contained in a cubic curve. We follow the approach of Green
and Tao in \cite{GT13}. The simple argument showing that $P \cup R$ is contained
in a cubic curve relies on Chasles theorem:

\begin{theorem}[Chasles Theorem]\label{theorem:Chasles}
Let $\ell_{1}, \ell_{2}, \ell_{3}$ be three distinct lines that
meet another family of three distinct lines $m_{1}, m_{2}, m_{3}$
at nine distinct intersection points.
Then every cubic curve passing through $8$ of these intersection points
must pass also through the $9$'th.
\end{theorem}

As a direct consequence of Theorem \ref{theorem:Chasles} we get the
following.

\begin{claim}\label{claim:1}
For $i\geq 0$,
Any cubic passing through eight of the points
$x_{n-1-i}, x_{n-2-i}, \ldots, x_{n-6-i}$
and $r_{5+2i}, r_{7+2i}, r_{9+2i}$ must pass through all nine of them.
\end{claim}

\begin{figure}[tb]
\begin{center}
\scalebox{0.6}{
\includegraphics{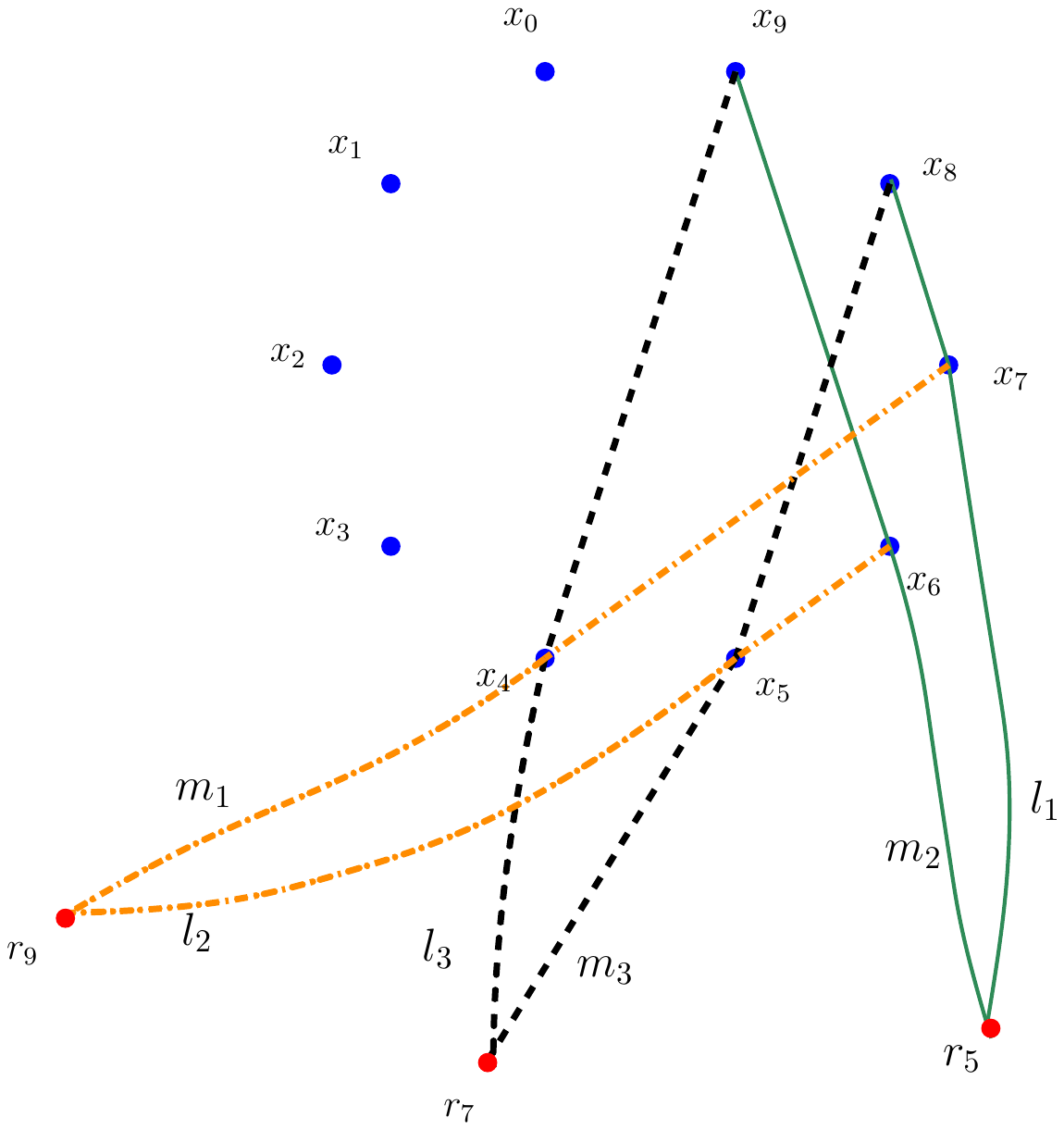}
}
\caption{Illustration for the proof of Claim~\ref{claim:1} with $n=10$ and $i=0$.}
\label{fig:cubic1}
\end{center}
\end{figure}

\noindent {\bf Proof.}
This is a direct application of Theorem \ref{theorem:Chasles}. Take
$\ell_{1}$ to be the line through $x_{n-2-i}, x_{n-3-i}$, and $r_{5+2i}$.
Take $\ell_{2}$ to be the line through $x_{n-4-i}, x_{n-5-i}$, and $r_{9+2i}$.
Take $\ell_{3}$ to be the line through $x_{n-1-i}, x_{n-6-i}$, and $r_{7+2i}$.
Take $m_{1}$ to be the line through $x_{n-3-i}, x_{n-6-i}$, and $r_{9+2i}$.
Take $m_{2}$ to be the line through $x_{n-1-i}, x_{n-4-i}$, and $r_{5+2i}$.
Take $m_{3}$ to be the line through $x_{n-2-i}, x_{n-5-i}$, and $r_{7+2i}$. (See Figure~\ref{fig:cubic1})

Notice that the nine intersection points of a line from
$\ell_{1}, \ell_{2}, \ell_{3}$ and a line from $m_{1}, m_{2}, m_{3}$
are $x_{n-1-i}, x_{n-2-i}, \ldots, x_{n-6-i}$,
and $r_{5+2i}, r_{7+2i}$, and $r_{9+2i}$.
The assertion of the claim follows now from
Theorem \ref{theorem:Chasles}.
\bbox

We use the easy and well known fact that for every nine points in the plane
there is a cubic curve passing through all of them.
Let $\Gamma$ be a cubic curve passing through
$x_{n-1}, x_{n-2}, \ldots, x_{n-7}$ and also through $r_{5}$ and $r_{7}$.
A repeated application of Claim \ref{claim:1} shows that
$\Gamma$ must pass through all the points $x_{0}, \ldots, x_{n-1}$
as well as through any point of the form $r_{5+2i}$. If $n$ is odd, then
this implies that all the points of $P \cup R$ are on $\Gamma$ and we are done.
If $n$ is even this only implies that all the points of $P$ and the points $r_{1}, r_{3}, \ldots, r_{n-1}$
are on $\Gamma$ and we need to show that also $r_{0}, r_{2}, \ldots, r_{n-2}$
lie on $\Gamma$. This will follow from the following claim.

\begin{claim}\label{claim:2}
Any cubic passing through eight of the points
$x_{n-1-i}, x_{n-2-i}, x_{n-3-i}, x_{n-5-i},
x_{n-6-i}, x_{n-7-i}$
and $r_{7+2i}, r_{8+2i}, r_{9+2i}$ must pass through all nine of them.
\end{claim}

\begin{figure}[tb]
\begin{center}
\scalebox{0.6}{
\includegraphics{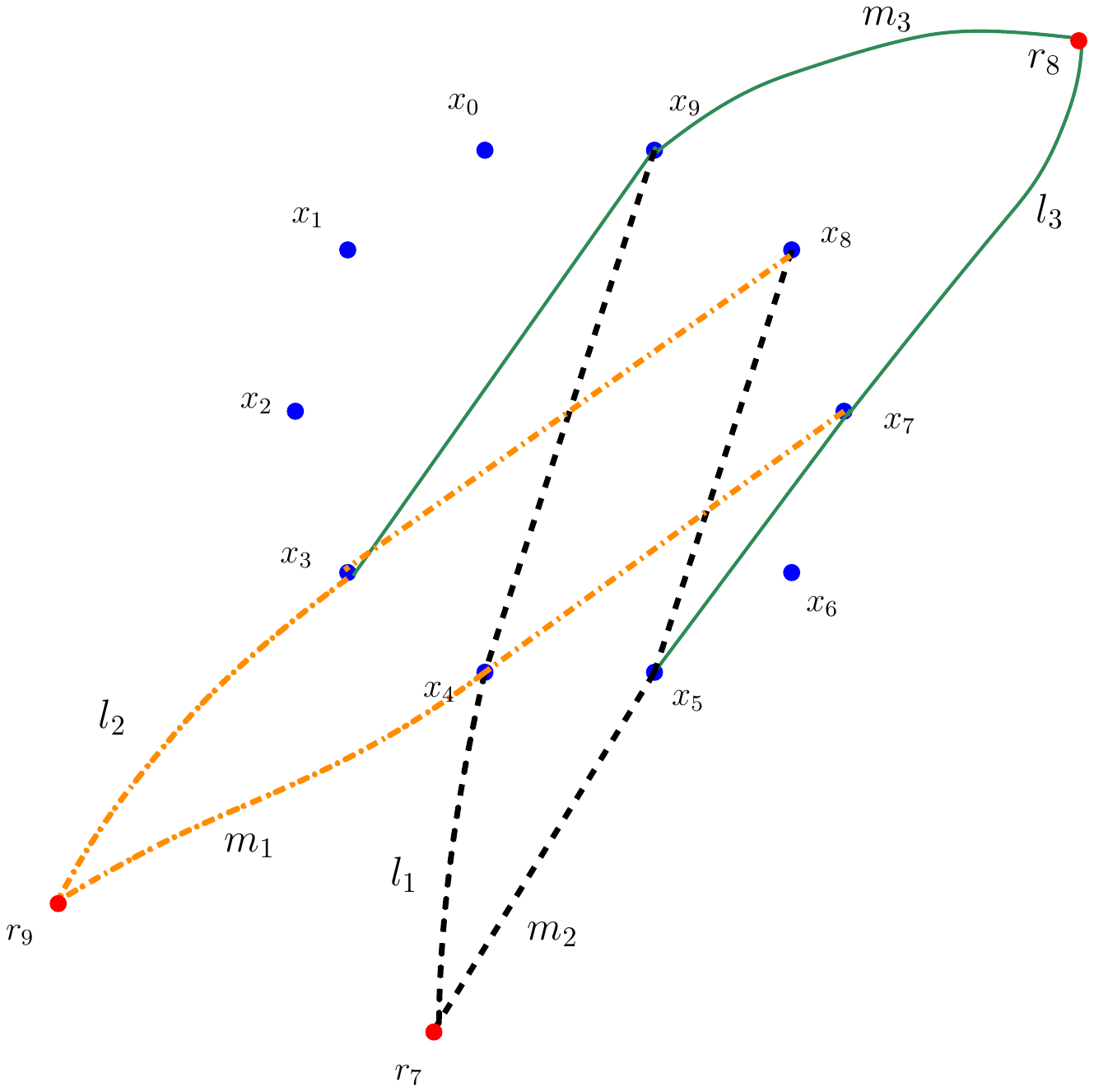}
}
\caption{Illustration for the proof of Claim~\ref{claim:2} with $n=10$ and $i=0$.}
\label{fig:cubic2}
\end{center}
\end{figure}

\noindent {\bf Proof.}
This is a direct application of Theorem \ref{theorem:Chasles}. Take
$\ell_{1}$ to be the the line through $x_{n-1-i}, x_{n-6-i}$, and $r_{7+2i}$.
Take $\ell_{2}$ to be the line through $x_{n-2-i}, x_{n-7-i}$, and $r_{9+2i}$.
Take $\ell_{3}$ to be the line through $x_{n-3-i}, x_{n-5-i}$, and $r_{8+2i}$.
Take $m_{1}$ to be the line through $x_{n-3-i}, x_{n-6-i}$, and $r_{9+2i}$.
Take $m_{2}$ to be the line through $x_{n-2-i}, x_{n-5-i}$, and $r_{7+2i}$.
Take $m_{3}$ to be the line through $x_{n-1-i}, x_{n-7-i}$, and $r_{8+2i}$. (See Figure~\ref{fig:cubic2})

Notice that the nine intersection points of a line from
$\ell_{1}, \ell_{2}, \ell_{3}$ and a line from $m_{1}, m_{2}, m_{3}$
are $x_{n-1-i}, x_{n-2-i}, x_{n-3-i}, x_{n-5-i}, x_{n-6-i}, x_{n-7-i}$,
and $r_{7+2i}, r_{8+2i}, r_{9+2i}$.
The assertion of the claim follows now from
Theorem \ref{theorem:Chasles}.
\bbox

\section{A bipartite version of Theorem \ref{theorem:main}}

In this section we consider a bipartite version of Theorem \ref{theorem:main}
in which $P$ is a set of $2n$ points that is the union of a set $B$ of
$n$ blue and a set $G$ of $n$ green points. The set $R$ is a set of $n$ red
points where we assume that the sets $B, G$, and $R$ are pairwise disjoint.
We also assume that the set
$P=B \cup G$ is in general position in the sense that no
three of its points are collinear. Assume that every line through a point in $B$
and a point in $G$ contains also a point from $R$. We would like to show that
$P \cup R$ lies on a cubic curve. We are able to show this under the
additional assumption that whenever $b \in B$, $g \in G$, and $r \in R$
are collinear, then $r$ does not lie between $b$ and $g$.

\begin{theorem}\label{theorem:bipartite}
Let $B, G$, and $R$ be three pairwise disjoint sets of points in the
plane. Assume that $B \cup G$ is in general position and
$|B|=|G|=|R|=n$. If every line through a point $b \in B$ and a point
$g \in G$ contains a point $r \in R$ that does not lie between $b$ and $g$,
then $B \cup G \cup R$ lies on some cubic curve.
\end{theorem}

\noindent {\bf Proof.}
Here too the proof will consist of two steps. We first show that $B \cup G$
must be in convex position in such a way that the points of $B$ and $G$
are arranged alternately on the boundary of the convex hull of $B \cup G$.
We will also be able to uniquely determine the structure of the
collinear triples of a point from $B$, a point from $G$, and a point from $R$.
In the second step we will use Chasles theorem to conclude that
$B \cup G \cup R$ lie on a cubic curve.

One easy but crucial observation is that any line through two points of
different colors in $B \cup G \cup R$ must pass through through a unique
point of the third color. This is because $B \cup G$ is in general position,
$|B|=|G|=|R|=n$,
and there is a point in $R$ on each one of the $n^2$ lines through a blue
point and a green point.

We split into two possible cases:

\noindent {\bf Case 1.} On the boundary of the convex hull of
$B \cup G$ there are only points of the same color, without loss of
generality let it be blue. We will show that this case is impossible.
Let $\ell$ be a line supporting the convex hull of $B \cup G$ in an edge
whose vertices are $b_{1}, b_{2} \in B$. Rotate $\ell$ about $b_{1}$
so that it separates $b_{2}$ from the rest of the points in $B \cup G$
until the first time it meets a point $g$ in $G$. There is a point $r \in R$
on the line through $b_{1}$ and $g$, and $g$ must lie between $b_{1}$ and $r$
on that line. To obtain a contradiction observe that the line through
$r$ and $b_{2}$ cannot contain a point $g'$ of $G$ in such a way that
$r$ does not lie between $b_{2}$ and $g'$ on that line.

\medskip

\noindent {\bf Case 2.} On the boundary of the convex hull of $B \cup G$
there are points of both colors blue and green.
We claim that along the boundary of the convex hull of $B \cup G$
blue and green points appear alternately.
To see this, notice that the contrary assumption is, without loss of generality,
that one can find three consecutive vertices along the convex hull of
$B \cup G$ such that the first one $b_{1}$ is blue while the next two
$g_{1}$ and $g_{2}$ are green.
Consider the line $\ell$ through $b_{1}$ and $g_{2}$. It must contain a point
$r \in R$ that lies outside the convex hull of $B \cup G$.
We get a contradiction as the line through $r$ and $g_{1}$ cannot contain
any point of $B$.

Next we show that the points of $B \cup G$ must be in convex position.
Denote by $k$ the number of blue points that is also the same as the number
of green points on the boundary of the convex hull of $B \cup G$.
Notice that every point of $R$ outside the convex hull of $B \cup G$
must be collinear with unique two edges of the convex
hull of $B \cup G$. Therefore, there are precisely $k$ points
of $R$ outside the convex hull of $B \cup G$.
Moreover, these $k$ points of $R$ are exactly those points of $R$
that pierce all the $k^2$ lines through two points of different
colors on the boundary of the convex hull of $B \cup G$.
Therefore, if $r \in R$ lies outside the convex hull of $B \cup G$, then
any line passing through $r$ and a point on the boundary of the
convex hull of $B \cup G$ must pass through another point of different color
also lying on the boundary of the convex hull of $B \cup G$.

\begin{figure}[tb]
\begin{center}
\scalebox{0.6}{
\includegraphics{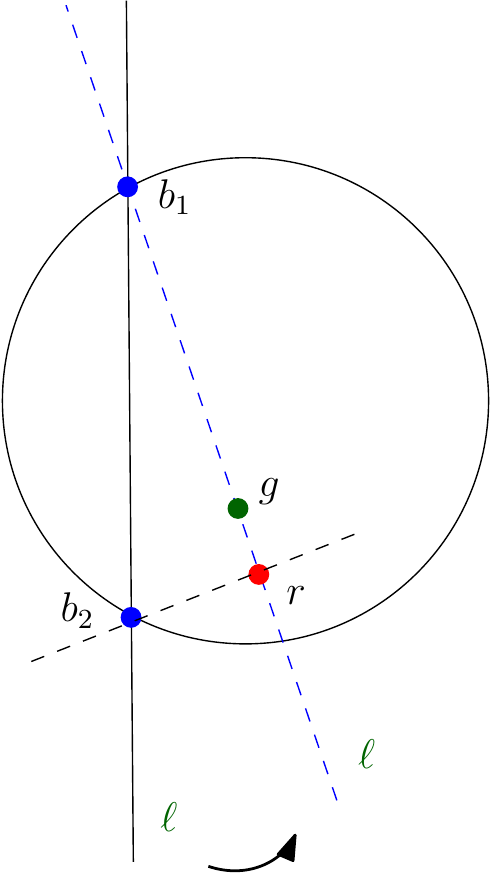}
}
\caption{Illustration for the proof of Theorem~\ref{theorem:bipartite}.}
\label{fig:empty2}
\end{center}
\end{figure}

In order to show that the points of $B \cup G$ are in convex position
it is enough to show that there is no point of $G$
inside the convex hull of $B \cup G$.
Let $b_{1}$ and $b_{2}$ be two blue points on the boundary of
the convex hull of $B \cup G$. Let $\ell$ be the line through $b_{1}$ an $b_{2}$. (See Figure~\ref{fig:empty2})
It is enough to show that as we rotate $\ell$ in the counterclockwise direction
around $b_{1}$, as long as $\ell$ does not contain an edge of the convex
hull of $B \cup G$, it never meets a point of $G$ inside the convex hull of
$B \cup G$ (then we can change the role of $b_{1}$ and $b_{2}$).
Assume to the contrary that $g \in G$ is the first such point let
$r$ be the point
of $R$ on the line through $b_{1}$ and $g$. The point $r$ cannot lie outside of
the convex hull of $B \cup G$ because the line through $r$ and $b_{1}$ does
not contain a green point on the boundary of the convex hull of $B \cup G$.
Now consider the line through $r$ and $b_{2}$. We get a contradiction
as this line cannot contain a green point $g'$ in such a way that $r$
does not lie between $b_{2}$ and $g'$ on that line.

Having shown that the points of $B \cup G$ lie in convex position
and alternately, let us denote the points of $B \cup G$ by
$x_{0}, x_{1}, \ldots, x_{2n-1}$ according to their cyclic order on the boundary
of the convex hull of $B \cup G$. We claim that for every odd number $k$
between $0$ and $2n$ there is a unique point $r$ (that we will denote
by $r_{2n-k}$) such that $r$ is collinear with all the pairs $x_{i}$ and $x_{j}$
such that $i+j=k$ (modulo $2n$).
Indeed, under the contrary assumption one can find $i$ and $j$ such that
$i+j=k$, $x_{i}$ and $x_{j}$ are collinear with $r \in R$ but at the same
time $r$ is not collinear with $x_{i+1}$ and $x_{j-1}$.
To get a contradiction rotate the line $\ell$ through $x_{i}$ and $x_{j}$
about $r$ such that it separates $x_{i}$ and $x_{j}$ from $x_{i+1}$ and $x_{j-1}$
until the first time it hits a point in $B \cup G$.
This point must be either $x_{i+1}$ or $x_{j-1}$ and no other point of
$B \cup G$ lies on that line, a contradiction.

Therefore, after denoting the points of $R$ by
$r_{1}, r_{3}, \ldots, r_{2n-1}$, as suggested above, we see that
$x_{i}, x_{j}$, and $r_{k}$ are collinear iff $i+j+k=0$ modulo $2n$.

Now comes the second part of the proof where we will use Chasles theorem
(in the form of Claim \ref{claim:1})
to show that $B \cup G \cup R$ lie on a cubic.

Let $\Gamma$ be a cubic curve passing through $x_{2n-1}, x_{2n-2}, \ldots, x_{2n-7}$
as well as through $r_{5}$ and $r_{7}$.
A repeated application of Claim \ref{claim:1} (with $2n$ in the role of $n$
in the statement of Claim \ref{claim:1}) shows that
all the points $x_{0}, x_{1}, \ldots, x_{2n-1}$ and
$r_{1}, r_{3}, r_{5}, \ldots, r_{2n-1}$ must lie on $\Gamma$.
\bbox

\bibliographystyle{alpha}
\bibliography{references}

\end{document}